\documentclass[11pt]{article}%
\usepackage{amssymb}
\newtheorem{theorem}{Theorem}[section]%
\newtheorem{lemma}[theorem]{Lemma}%
\newtheorem{cor}[theorem]{Corollary}%
\newtheorem{exam}[theorem]{Example}%
\newtheorem{rem}{Remark}%
\def\a{\alpha}

 \def\O{\Omega}

 \def\og{\overline G} 
 \def\oq{\overline Q} 
  
 \def\op{\overline P}

\def\lg{\langle} \def\rg{\rangle}

\def\di{\bigm|}

\def\f{\noindent}

\def\qed{\hfill $\Box$}

\def\demo{{\bf Proof}\hskip10pt}
\def\A{\mathcal{A}}

\def\C{\hbox{\rm C}}

\def\Aut{\hbox{\rm Aut}}

\def\Syl{\hbox{\rm Syl}}

\def\mod{\hbox{\rm mod }}

\begin{document}
\baselineskip=16pt
\title{Finite groups with few conjugate classes of minimal non-abelian subgroups
\thanks{This work was supported by  NSFC(No. 12271318 \& 12371022) and FRPS(No. 202203021221126).}}
\author{\\ Haipeng Qu and Junqiang Zhang\thanks{Corresponding author. e-mail: junqiangchang@163.com}\\
Department of Mathematics , Shanxi Normal University\\
Taiyuan Shanxi 030031 China\\
\\}
\date{}
\maketitle

\begin{abstract}
Let $G$ be a finite non-abelian group and $\kappa_1(G)$ the number of conjugate classes of minimal non-abelian subgroups of $G$. The structure of $G$ with $\kappa_1(G)=1$ is determined. In the case of $G$ being the $p$-groups, the structure of $G$ with $\kappa_1(G)\leqslant p$ is also determined.

\medskip
\vspace{0.2cm}

 \f {\bf Key Words:} minimal non-abelian subgroups; conjugate classes of subgroups; Frobenius groups; Fermat primes;  Mersenne primes.

\medskip
\f {\it \bf AMS Subject Classifications:} 20D15, 20D20.
\end{abstract}

\section{Introduction}

All groups considered in this paper are finite. A finite non-abelian group is said to be minimal non-abelian if its proper subgroups are abelian.
Obviously, every finite non-abelian group contains at least one minimal non-abelian subgroup. Therefore, it is natural to study
the structure of finite groups by minimal non-abelian subgroups.
Assume $G$ is a finite non-abelian group and $\A_1(G)$ the set consisting of all minimal non-abelian subgroups of $G$.
Let $G$ act conjugately on $\A_1(G)$ and  $\kappa_1(G)$ denote the number of orbits of $G$ on $\A_1(G)$. In other words, $\kappa_1(G)$ is the number of conjugate classes of minimal non-abelian subgroups in $G$.

The aim in this paper is to determine the structure of $G$ with $\kappa_1(G)=1$. We prove (see Theorem \ref{MT1ml}) that if $G$ has a non-abelian Sylow subgroup,
then $\kappa_1(G)=1$ if and only if $\alpha_1(G)=1$, where $\alpha_1(G)$ is the number of minimal non-abelian subgroups of $G$.
Fortunately, such groups with $\alpha_1(G)=1$ have been determined by Berkovich, which are the groups (a) or (b) in \cite[Theorem 4.1]{Be1}.

Assume $\alpha_1(G)>1$ and all Sylow subgroups of $G$ are abelian. Under such assumption, such groups $G$ with $\kappa_1(G)=1$ are determined, see [Theorem \ref{MT2-dl1}]. Furthermore, we prove (see Theorem \ref{MT2-dl2}) that if $\kappa_1(G)=1$, then $G/Z(G)$ is a Frobenius group, where the Frobenius kernel is a homocyclic $p$-subgroup, the Frobenius complement is a cyclic $q$-subgroup, where $p, q$ are distinct primes. Moreover, if $G$ is $2$-transitive on $\A_1(G)$, then $q=2$ and $p$ is a Fermat prime, or $p=2$ and $q$ is a Mersenne prime.

In the case of $G$ being non-abelian $p$-group, if $\kappa_1(G)=1$, then, as a direct consequence of \cite[Proposition 10.28]{BE1}, $G$ itself is minimal non-abelian. If  $\kappa_1(G)>1$, then, it is easily follows by \cite[Theorem 17]{BE} that $\kappa_1(G)\geqslant p$. We prove (see Theorem \ref{Main}) that $\kappa_1(G)=p$ if and only if $G$ has an abelian subgroup of index $p$ and all non-abelian subgroups of $G$ are generated by two elements. Such $p$-groups were classified by Xu et al., see \cite[Theorem 3.12 and 3.13]{XAZ}.

\medskip
The terminology and notations are standard, as in \cite{BE1}.



%
%
%
%


\section{Finite groups $G$  with $\kappa_1(G)=1$}

The structure of minimal non-abelian groups has been determined by Miller and Moreno in \cite{Miller2}.
They proved the following theorem, but they didn't list explicitly it. For convenience we list it as follows.

\begin{theorem}\label{MM} Assume $G$ is a non-nilpotent minimal non-abelian group. Then 
$G=P\rtimes Q$, where $Q\in {\rm Syl}_q(G)$ is cyclic, $G'=P
\in {\rm Syl}_p(G)$ is an
elementary abelian subgroup, $p$ and $q$ are distinct primes.
\end{theorem}

From their result we observe
that the commutator subgroup of any minimal non-abelian group is a $p$-subgroup for some prime number $p$. Let $\pi_1(G)$ be the set consisting of these prime numbers. In other words,
$$\pi_1(G)=\{ p \mid \exists H\in\A_1(G)\ \mbox{such that}\ H' \ \mbox{is a $p$-subgroup}\}.$$

\begin{theorem}\label{pi1}
Assume $G$ is a finite group. Then $G=K\rtimes A$, where $K$ is a $\pi_1(G)$-subgroup and $A$ is an abelian $\pi_1(G)'$-subgroup.
In particular,

{\rm (1)} if $\kappa_1(G)=1$, then $G=P\rtimes A$, where $P$ is a $p$-subgroup and $A$ is an abelian $p'$-subgroup of $G$.

{\rm (2)} if $\kappa_1(G)\leqslant 2$, then $G$ is solvable. If $\kappa_1(G)\geqslant 3$, then $G$ is not necessarily solvable. For example, $\kappa_1(A_5)=3$.
\end{theorem}

\demo Obviously, $\pi_1(G)'$-subgroups of $G$ are abelian.
Let $p\not\in \pi_1(G)$ and $P$ is a Sylow $p$-subgroup of $G$. Then $P$ is abelian.

We assert that $N_G(P)=C_G(P)$.
If not, there exists $x\in N_G(P)\setminus C_G(P)$. Since $P$ is abelian, $P\leq C_G(P)$.
Assume that $x$ is a $p'$-element without loss of generality.
Then $G$ has a non-abelian subgroup $L=P\rtimes \lg x\rg$.
It follows that $H'\leq L'\leq P$ for any $H\in\A_1(L)$. Thus $H'$ is a $p$-subgroup of $G$.
This contradicts $p\not\in\pi_1(G)$. Hence $N_G(P)=C_G(P)$.

By Burnside's normal $p$-complement theorem, $G$ has a normal $p'$-subgroup $K_p$ such that $G=K_p\rtimes P$.
Let $K=\bigcap\limits_{p\not\in\pi_1(G)}K_p$. Then $K$ is a normal $\pi_1(G)$-Hall subgroup of $G$.
By Schur-Zassenhaus Theorem, $G$ has a subgroup $A$ such that $G=K\rtimes A$.
It follows that $A$ is a $\pi_1(G)'$-Hall subgroup and so $A$ is abelian.\qed

\begin{rem}
{\rm (1)} If we replace $\kappa_1(G)=1$ in Theorem $\ref{pi1}(1)$ by the condition ``minimal non-abelian subgroups of $G$ are isomorphic", then the conclusion still holds.

{\rm (2)} If minimal non-abelian subgroups of $G$ are of same order, then $|\pi_1(G)|\leqslant 2$ in Theorem $\ref{pi1}$.
\end{rem}





\begin{lemma}\label{zhij}
Assume a finite group $G=G_1\times G_2$, where $(|G_1|,|G_2|)=1$. Then
\begin{itemize}
\item[{\rm (1)}] $\alpha_1(G)=\alpha_1(G_1)+\alpha_1(G_2);$

\item[{\rm (2)}]  $\kappa_1(G)=\kappa_1(G_1)+\kappa_1(G_2)$.
\end{itemize}
\end{lemma}

\demo For any $H\leq G$, there exist $H_1\leq G_1$ and $H_2\leq G_2$ such that $H=H_1\times H_2$. Thus any minimal non-abelian subgroup of $G$ is in $G_1$ or $G_2$ and so the results hold. \qed

\medskip
The following Lemma \ref{geqp} can be easily followed by \cite[Proposition 10.28]{BE1} and  by \cite[Theorem 17]{BE}, respectively. We list explicitly them without proof.

\begin{lemma}\label{geqp}
Let $G$ be a finite $p$-group.
\begin{itemize}
\item[{\rm (1)}]  If $\kappa_1(G)=1$, then $G$ is minimal non-abelian;

\item[{\rm (2)}] If $\kappa_1(G)>1$, then $\kappa_1(G)\geqslant p$.
\end{itemize}
\end{lemma}

\begin{lemma}{\rm (\cite{Miller})}\label{Miller}
Let $G$ be a non-abelian $p$-group. Then the number of abelian subgroups of index $p$ in $G$ is $0,1$ or
$p+1$.
\end{lemma}

\begin{lemma}\label{t3(G)}
Let $G$ be a $p$-group of order $\geqslant p^4$. Then the number of non-abelian subgroups of order $p^3$ is derived by $p$.
\end{lemma}

\demo Assume $|G|=p^n$ and $t_i(G)$ denotes the number of non-abelian subgroups of order $p^i$ of $G$.  We use induction on $n$.
If $n=4$, then $G$ has an abelian subgroup of order $p^3$.
By Lemma \ref{Miller}, the number of abelian subgroups of order $p^3$ in $G$ is $1$ or $1+p$.
 On the other hand, the number of subgroups of order $p^3$ of $G$ is $1+p+\cdots+p^{d(G)-1}$. Thus $t_3(G)\equiv 0(\mod p)$. The result holds for $n=4$.

Assume $n>4$ and the result holds for any maximal subgroup $M$ of $G$. That is, $t_{3}(M)\equiv 0~(\mod p)$. By the enumeration principle of P. Hall (see \cite[Theorem 5.2]{BE1}),
$t_{3}(G)\equiv{\sum \limits_{M\lessdot G}{t_{3}(M)}}\equiv 0~(\mod p)$, where $M\lessdot G$ denotes $M$ is a maximal subgroup of $G$.  The result holds. \qed

\begin{lemma}{\rm (\cite[Lemma 2.2]{XAZ})}\label{minimal non-abelian equivalent conditions}
Let $G$ be a finite $p$-group. Then the following conditions are
equivalent:
\begin{itemize}
\item[{\rm (1)}]  $G$ is minimal non-abelian;

\item[{\rm (2)}] $d(G)=2$ and $|G'|=p$;

\item[{\rm (3)}]  $d(G)=2$ and $\Phi(G)=Z(G)$.
\end{itemize}
\end{lemma}

\begin{lemma}\label{p3} Let $G$ be a finite minimal non-abelian $p$-group.
If $\exp(G)\leqslant p^2$ and $\Phi(G)$ is cyclic, then $|G|=p^3$ .
\end{lemma}

\demo By Lemma \ref{minimal non-abelian equivalent conditions}, $d(G)=2$ and $|G'|=p$.
Let $G=\lg a,b\rg$. Then $\Phi(G)=\lg a^p,b^p,[a,b]\rg$.
It follows by $\exp(G)\leqslant p^2$ that $\exp(\Phi(G))=p$. Since $\Phi(G)$ is cyclic, $\Phi(G)\cong\C_p$.
Since $d(G)=2$, $|G/\Phi(G)|=p^2$. So  $|G|=p^3$. \qed

\begin{lemma}\label{A1d2} Let $G$ be a finite abelian $p$-group and $1<H\le G$. Then there exists $N\le G$ such that $G/N$ is cyclic and
$|H:H\cap N|=p$.
\end{lemma}

\demo If $G$ is cyclic, then the maximal subgroup of $H$ is the required subgroup.  Assume $G=\lg a_1\rg\times \lg a_2\rg\times \cdots \times \lg a_n\rg$, where $n\geqslant 2$.
Let $N_i=\lg a_1,a_2,\cdots,a_{i-1},a_{i+1},\cdots, a_n\rg$. Then $\bigcap\limits_{i=1}^n N_i=1$.
So there exists $i\in\{1,2,\cdots,n\}$ such that $H\nleq N_i$. Clearly, $G/N_i$ is cyclic and
$HN_i/N_i\ne 1$. Let $N/N_i$ be a maximal subgroup of $HN_i/N_i$. Then $N$ is the required subgroup. \qed

\begin{theorem} \label{MT1ml}
Assume  $G$ is a finite group with a non-abelian Sylow subgroup.
Then $\kappa_1(G)=1$ if and only if $\alpha_1(G)=1$, where $\alpha_1(G)$ is the number of minimal non-abelian subgroups in $G$.
\end{theorem}

\demo $\Longleftarrow:$ Obviously.

$\Longrightarrow:$
Let $P$ be a non-abelian Sylow $p$-subgroup of $G$. By Theorem \ref{pi1}(1), $G=P\rtimes A$, where $A$ is an abelian $p'$-subgroup.
If $G$ is nilpotent, then $G=P\times A$.
By Lemma \ref{zhij}(2),  $\kappa_1(G)=\kappa_1(P)+\kappa_1(A)=\kappa_1(P).$
Thus $\kappa_1(G)=1$ if and only if $\kappa_1(P)=1$. By Lemma \ref{geqp}(1), $P$ is minimal non-abelian.
By Lemma \ref{zhij}(1), $\alpha_1(G)=\alpha_1(P)+\alpha_1(A)=\alpha_1(P)=1.$

\medskip
Assume $G$ is non-nilpotent.  We prove that $\alpha_1(G)=1$ by following six steps.

\medskip
{\bf Step 1.} All minimal non-abelian subgroups of $G$ are contained in $P$.

\smallskip
Since $P$ is not abelian, $P$ has a minimal non-abelian subgroup $H$.
Since $\kappa_1(G)=1$ and $P\unlhd G$, all minimal non-abelian subgroups of $G$ are in $P$.

\medskip
{\bf Step 2.} $N$ is abelian and $[N,A]=1$ for any $G$-invariant subgroup $N< P$.
\smallskip

If $N$ is not abelian, then $N$ has a minimal non-abelian subgroup $H_1$.
Since $N<P$, by \cite[Proposition 10.28]{BE1}, $P$ has a minimal non-abelian subgroup $H_2$ such that $H_2\nleq N$.
It follows from $N\unlhd G$ that $H_2$ is not conjugate to $H_1$. This contradicts that $\kappa_1(G)=1$. Thus $N$ is abelian.
If $[N,A]\neq 1$, then there is a minimal non-abelian subgroup $H\leq N\rtimes A$ such that $H$ is not a $p$-subgroup, which contradicts the result of Step 1. Thus $[N,A]=1$.

\medskip
{\bf Step 3.} $A^G=G$.
\smallskip

If $A^G<G$, then $P\cap A^G< P$. Let $N_1=P\cap A^G$. Then $N_1\unlhd G$. By Step 2, $N_1$ is abelian and $[N_1,A]=1$.
Let $\overline{G}=G/N_1$. Notice that $$[G,A]\leq G'\cap A^G\leq P\cap A^G=N_1.$$ We get $\overline{G}=\overline{P}\times \overline{A}$, which is nilpotent.
Thus $K_n(G)\leq N_1$ for some positive integer $n$. Notice that $N_1\leq P$. We get $[N_1,mP]=1$ for some positive integer $m$.
Since $[N_1,A]=1$, $[N_1,mG]=1$.  So $K_{m+n}(G)=1$. This contradicts $G$ is non-nilpotent. Thus $A^G=G$.

\medskip
{\bf Step 4.}  $\Phi(P)\leq Z(G)$, $\exp(P')=p$, $\exp(P)=p$ when $p>2$, $\exp(P)=4$ when $p=2$.
\smallskip

For any $G$-invariant subgroup $N< P$, $N$ is abelian and $[N,A]=1$ by Step 2. Thus $[N,G]=[N,A^G]=1$ by Step 3 and so $N\leq Z(G)$.

Notice that $\Phi(P)$ is a $G$-invariant subgroup and $\Phi(P)< P$. We get $\Phi(P)\leq Z(G)$. It follows that $\exp(P')=p$.
Since $G$ is non-nilpotent, $A$ acts non-trivially on $P$.
By \cite[Kapitel VI, Satz 5.12]{Hup}, $P=\Omega_1(P)$ when $p>2$ or $P=\Omega_2(P)$ when $p=2$.
It follows that $\exp(P)=p$ when $p>2$ or $\exp(P)=4$ when $p=2$.

\medskip
{\bf Step 5.}  
$\A_1(P/N)$ is an orbit of $G/N$ if $P'\nleq N\le \Phi(P)$.
\smallskip

Let $\overline{P}=P/N$ and $\overline{G}=G/N$.
Clearly, $\A_1(\op)$ is $\og$-invariant. It suffices to show $\og$ is transitive on $\A_1(\op)$.
For any $H\in\A_1(\overline{P})$,
there exists $K=\lg a,b\rg\leq P$ such that $H=KN/N$.
Since $\Phi(P)\leq Z(G)$, $\Phi(K)\leq Z(K)$. Notice that $K$ is non-abelian. We get $Z(K)\leq \Phi(K)$ and so $Z(K)=\Phi(K)$.
By Lemma \ref{minimal non-abelian equivalent conditions}, $K\in\A_1(P)$.
This means that every minimal non-abelian subgroup of $\overline{P}$ is an image of a minimal non-abelian subgroup of $P$.
It follows by $\kappa_1(G)=1$ that $\og$ acts transitively on $\A_1(\op)$.

\medskip
{\bf Step 6.}  $\alpha_1(G)=1$.
\smallskip

By Step 1, it suffices to show that $P\in\A_1(G)$. By Lemma \ref{minimal non-abelian equivalent conditions}, it suffices to show that $\Phi(P)=Z(P)$ and $d(P)=2$.
By Step 4, $\Phi(P)\leq Z(G)$ and so $\Phi(P)\leq Z(P)$. If $d(P)=2$, then $Z(P)\leq \Phi(P)$ and so $\Phi(G)=Z(G)$. Thus it is enough to show $d(P)=2$.

By Lemma~\ref{A1d2}, there exists $N<\Phi(P)$ such that $\Phi(P)/N$ is cyclic and $|P':P'\cap N|=p$.
Let $\og=G/N$. Then $\Phi(\op)$ is cyclic and $|\op'|=p$. For any $L\in\A_1(\op)$, we have $L'=\op'$ and so $L\unlhd \op$.
Thus $\op\leq N_{\og} (L)$.
By Step 5, $\A_1(\op)$ is an orbit of $\og$ and so $$\alpha_1(\op)=|\og:N_{\og} (L)|\not\equiv 0(\mod p).$$

On the other hand, by Lemma \ref{p3}, $|L|=p^3$. If $|\op|\geqslant p^4$, then $\alpha_1(\op)\equiv 0(\mod p)$ by Lemma \ref{t3(G)}.
This is a contradiction. Thus $|\op|=p^3$ and so $d(P)=d(\op)=2$. \qed


\medskip
Berkovich determined the finite groups $G$ with $\alpha_1(G)=1$, see \cite[Theorem 4.1]{Be1}.
If $G$ has a non-abelian Sylow $p$-subgroup, then we know from \cite[Theorem 4.1]{Be1} that $\kappa_1(G)=1$ if and only if $G$ is the groups (a) or (b) in \cite[Theorem 4.1]{Be1}.

\medskip
It is needed to be pointed out that if all Sylow subgroups of finite groups $G$ are abelian, then $\kappa_1(G)=1$ is not necessarily equivalent to $\alpha_1(G)=1$.


\begin{exam}  Let $G=\lg a,b,c\di a^3=b^3=c^8=1, [a,b]=1, a^c=b, b^c=ab\rg$. Then $|G|=2^3\cdot 3^2$.
It is easy to get that $H\in\A_1(G)$ if and only if $H\in\A_1(\lg a,b,c^4\rg)$ and $\alpha_1(G)=12$.
On the other hand, $\lg a,c^4\rg\in\A_1(G)$ and $N_G(\lg a,c^4\rg)=\lg a,c^4\rg$.
It follows that $|G:N_G(\lg a,c^4\rg)|=12$ and so all minimal non-abelian subgroups of $G$ are conjugate each other. Thus $\kappa_1(G)=1$.
\end{exam}

In following, we determine the structure of finite groups $G$ whose Sylow subgroups are abelian and  $\kappa_1(G)=1$.


\begin{lemma} {\rm \cite[Chap 5, Theorem 2.3]{Gost}}\label{MT2-yl1}
Let $A$ be a $p'$-group of automorphism group of an abelian $p$-group $P$. Then
$$P=[P,A]\times C_P(A).$$
\end{lemma}

\begin{lemma} \label{G'A}
Let $G=G'\rtimes A$, where $G'$ is an abelian $p$-subgroup, $A$ is a $p'$-subgroup. Then $G'=[G',A]$ and $C_{G'}(A)=1$.
\end{lemma}

\demo By Lemma \ref{MT2-yl1}, $G'=[G',A]\times C_{G'}(A)$. Since $G'$ and $A$ are abelian, $G/[G',A]$ is abelian.
It follows that $G'=[G',A]$ and so $C_{G'}(A)=1$. \qed

\begin{lemma} \label{A1(B)}
Let $G=B\times S$, where $B=G'\rtimes Q$, $G'$ is an abelian $p$-subgroup, $Q$ is an abelian $q$-subgroup and $S$ is an abelian $q'$-subgroup, where $p, q$ are distinct primes.
Then  $\A_1(G)=\A_1(B)$.
\end{lemma}

\demo
For any $K\in\A_1(G)$, $K'\leq G'$.
By Theorem \ref{MM}, $K=K'\rtimes\lg b\rg$ for some $b\in Q^g$ and $g\in G$.
Since $B\unlhd G$, $Q^g\leq B$.
It follows that $K\leq B$ and so $\A_1(G)=\A_1(B)$. \qed

\begin{lemma} \label{CG(a)=1}
Let $G=P\rtimes Q$, where $P$ is an abelian $p$-subgroup and $Q$ is an  abelian $q$-subgroup, $p, q$ are distinct primes.
If $\Omega_1(P)$ is minimal normal, then $C_{P}(a)=1$ for any $a\in Q\setminus Z(G)$.
\end{lemma}

\demo Assume there exists $a\in Q\setminus Z(G)$ such that $C_{P}(a)\neq 1$. Then, by Lemma \ref{MT2-yl1}, $P=[P,\lg a\rg]\times C_{P}(\lg a\rg)$ and so $[P,\lg a\rg]\cap\Omega_1(P)<\Omega_1(P)$.
For any $b\in Q$, we have
$$[P,\lg a\rg]^b=[P^b,\lg a\rg^b]=[P,\lg a\rg].$$
Thus $Q\leq N_G([P,\lg a\rg])$. Noting $[P,\lg a\rg]\leq P$ and $P$ is abelian, we get $P\leq N_G([P,\lg a\rg])$ and so $[P,\lg a\rg]\unlhd G$.
Since $\Omega_1(G')$ is minimal normal, $[P,\lg a\rg]\cap\Omega_1(P)=1$. Thus $[P,\lg a\rg]=1$. It follows that $a\in Z(G)$. A contradiction.
\qed

\begin{theorem} \label{MT2-dl1}
Assume $G$ is a finite group whose Sylow subgroups are abelian.
If $\kappa_1(G)=1$, then $G=(G'\rtimes Q)\times S$, where $G'$ is a homocyclic $p$-subgroup, $\Omega_1(G')$ is a minimal normal subgroup of $G$,
$Q$ is a cyclic $q$-subgroup and $S$ is an abelian $q'$-subgroup, where $p, q$ are distinct primes.
\end{theorem}

\demo By Theorem \ref{pi1}(1), $G=P\rtimes A$, where $P$ is an abelian $p$-subgroup and $A$ is an abelian $p'$-subgroup.

\medskip
{\bf Step 1.} $G=(G'\rtimes Q)\times S$, where $G'$ is an abelian $p$-subgroup, $Q$ is an abelian $q$-subgroup, $S$ is an abelian $q'$-subgroup.

\smallskip
By Lemma \ref{MT2-yl1}, $P=[P,A]\times C_P(A)$.
It follows that $G=H\times C_P(A)$, where $H=[P,A]\rtimes A$.
Since all Sylow subgroups of $G$ are abelian, the order of a minimal non-abelian subgroup may be assumed by $p^aq^b$.
Thus $x\in C_A(P)$ for any $\{p,q\}'$-element $x$ and so there is an abelian $\{p,q\}'$-subgroup $A_1$ and an abelian $q$-subgroup $Q$ such that
$H=[P,Q]\rtimes Q\times A_1$. It follows that $G=[P,Q]\rtimes Q\times A_1 \times C_P(A)$.
Obviously, $G'=[P,Q]\leq P$. Thus $G'$ is an abelian $p$-subgroup. Let $S=A_1 \times C_P(A)$. Then $G=(G'\rtimes Q)\times S$.

\medskip
{\bf Step 2.} $\Omega_1(G')$ is a minimal normal subgroup of $G$.

\smallskip
If not, there exists $N< \Omega_1(G')$ such that $1<N \unlhd G$.
Since $G'$ is an abelian $p$-subgroup, by Maschke's theorem on complete reducibility, we get $\Omega_1(G')=N\times L$, where $L\unlhd G$.
By Lemma \ref{G'A}, $C_{G'}(A)=1$. There are $x,y\in A$ such that $[x,N]\neq 1$ and $[y,L]\neq 1$.
Thus there are $H_1\in\A_1(\lg x,N\rg)$ and $H_2\in\A_1(\lg y,L\rg)$ such that $H_1'\leq N$ and $H_2'\leq L$.
This implies that $H_1$ is not conjugate to $H_2$, which contradicts $\kappa_1(G)=1$.

\medskip
{\bf Step 3.} $G'$ is homocyclic.

\smallskip
Let $\exp(G')=p^e$. Then $1<\mho_{e-1}(G')\leq\Omega_1(G')$.
Since $\mho_{e-1}(G')$ is characteristic in $G'$, $\mho_{e-1}(G')\unlhd G$.
It follows by Step 2 that $\mho_{e-1}(G')=\Omega_1(G')$ and so $G'$ is homocyclic.

\medskip
{\bf Step 4.} $Q$ is cyclic.

\smallskip
Take $H_1\in\A_1(G)$, where $H_1=P_1\rtimes\lg x\rg$ for some $P_1\leq G'$ and $x\in Q$.
If $Q$ is not cyclic, then there is $y\in Q\setminus\lg x\rg$ such that $o(y)=q$. Consider $\og=G/G'$. Then
$\oq=QG'/G'\cong Q$. So $\lg \bar{x}\rg\cap \lg \bar{y}\rg=1$.
It follows that $(G'\rtimes\lg x\rg)\cap\lg y\rg=1$.

If $y\in Z(G)$, letting $H_2=P_1\rtimes\lg xy\rg$, then $H_2\in\A_1(G)$.
Since $\kappa_1(G)=1$, there exists $x_1\in H_1$ and $g\in G$ such that $x_1^g=xy$.
Then $$x^{-1}x_1[x_1,g]=y\in (G'\rtimes\lg x\rg)\cap \lg y\rg=1.$$ It follows that $y=1$. This contradicts $y\in Q\setminus\lg x\rg$.

If  $y\not\in Z(G)$, there exists $P_2\leq G'$ such that $P_2\rtimes\lg y\rg\in\A_1(G)$.
Since $\kappa_1(G)=1$, there exist $g_2\in G$ and $x_1\in H_1$ such that $x_1^{g_2}=y$.
It follows that $$x_1[x_1,g_2]=y\in (G'\rtimes\lg x\rg)\cap\lg y\rg=1.$$
It follows that $y=1$. A contradiction again. Thus $Q$ is cyclic. \qed

\medskip
Under the assumption as that of Theorem \ref{MT2-dl1}, we will prove that $G/Z(G)$ is a Frobenius group. Moreover, if $G$ is $2$-transitive on $\A_1(G)$, then the odd primes $p,q$ involved in a Frobenius group are exactly a Fermat prime or a Mersenne prime. It should be mentioned that the result of Catalan's conjecture is used. Catalan conjectured that the equation $x^u-y^v=1$ has no other solution in positive integers except $3^2-2^3=1$ for
$u>1, v>1$. Mih$\breve{a}$ilescu proved the conjecture is true, see \cite[Theorem 5]{Miha}.





\begin{theorem} \label{MT2-dl2}
Assume $G$ is a finite group whose Sylow subgroups are abelian. If $\kappa_1(G)=1$, then

{\rm (1)} $\kappa_1(G/Z(G))=1$ and $G/Z(G)$ is a Frobenius group,
where the Frobenius kernel is a homocyclic $p$-subgroup, the Frobenius complement is a cyclic $q$-subgroup.

{\rm (2)} if $G$ is $2$-transitive on $\A_1(G)$, then $q=2$ and $p$ is a Fermat prime, or $p=2$ and $q$ is a Mersenne prime.
\end{theorem}

\demo By Theorem \ref{MT2-dl1}, $G=(G'\rtimes Q)\times S$,
where $S\leq Z(G)$. Let $B=G'\rtimes Q$. Then $\A_1(G)=\A_1(B)$ by Lemma \ref{A1(B)}. So the transitivity of $G$ on $\A_1(G)$ is reduced to the transitivity of $B$ on $\A_1(G)$. Without loss of generality assume $G=G'\rtimes Q$.
It follows that $G'=[G',Q]$. By Lemma \ref{MT2-yl1}, $C_{G'}(Q)=1$. Thus $Z(G)=C_Q(G')$ and $G'\cap Z(G)=1$.

\medskip
We prove (1) in three steps.

\medskip
{\bf Step 1.}  $Z(G/Z(G))=1$.

\smallskip
Let $\overline{G}=G/Z(G)$ and $\overline{Q}\cong Q/Z(G)$. Then $\overline{G}\cong G'\rtimes \overline{Q}$.
Since $G'\cap Z(G)=1$,
$$[a,b]=1\ \mbox{if and only if}\ [\bar{a},\bar{b}]=\bar{1}\ \mbox{ for\ any}\ a,b\in G.$$
Thus for any $\bar{z}\in Z(\overline{G})$,
we have $z\in Z(G)$.
It follows that $Z(G/Z(G))=1$.

\medskip
{\bf Step 2.} $\kappa_1(G/Z(G))=1$.

\smallskip
Let $K\in \A_1(\overline{G})$. Then $K=K'\rtimes\lg \bar{x}\rg$ for some $\bar{x}\in \overline{Q}^{\overline{g}}$ by Theorem \ref{MM}. Thus $K'=K\cap \overline{G}'$.
Take $D\leq G'$ such that $DZ(G)/Z(G)=K'$. Since $G'\cap Z(G)=1$, $D\cap Z(G)=1$.
It follows that $D\cong D/(D\cap Z(G))\cong K'$.
Let $H=\lg D, x\rg$. Notice that the action of $x$ on $D$ is the same as the action of $\bar{x}$  on $\overline{K}'$.
We get $H=D\rtimes \lg x\rg\in\A_1(G)$.
That means every minimal non-abelian subgroup of $\overline{G}$ is an image of a minimal non-abelian subgroup of $G$.
It follows by $\kappa_1(G)=1$ that $\kappa_1(G/Z(G))=1$.

\medskip
{\bf Step 3.} $G/Z(G)$ is a Frobenius group.

\smallskip
By Step 1 and Step 2, we may assume $Z(G)=1$ and $G=G'\rtimes Q$. Then $C_{G'}(Q)=1$ and $\Omega_1(G')$ is a minimal normal subgroup of $G$.
By Lemma \ref{CG(a)=1}, $C_{G'}(a)=1$ for any $a\in Q\setminus \{1\}$. It follows that $G$ is a Frobenius group, where the Frobenius kernel $G'$ is a homocyclic $p$-subgroup, the Frobenius complement $Q$ is a cyclic $q$-subgroup.

\medskip
We prove (2) in two steps. Now, $|\A_1(G)|=\alpha_1(G)\geqslant 2$. So $H\ntrianglelefteq G$ for all $H\in \A_1(G)$.

\medskip
{\bf Step 4.} $\exp(G')=p^2$ and $N_G(H)\cap G'=H'=\Omega_1(G')<G'$ for each $H\in \A_1(G)$.

\smallskip
Without loss of generality, assume $H=H'\rtimes\lg x\rg$, where $x\in Q$. Then $H'\leq \Omega_1(G')$.
Let $N=N_G(H)$ and $P\in{\rm Syl}_p(N)$.
Notice that $G'\in{\rm Syl}_p(G)$ and $H'\in{\rm Syl}_p(H)$. We get $P=N\cap G'$ and $H'=H\cap G'$.
Thus $H'\leq P$. On the other hand, since $H'\leq P=N\cap G'$, $1\neq [x,P]\leq P$.
By Lemma \ref{MT2-yl1}, $P=[P,\lg x\rg]\times C_P(\lg x\rg)$.
Since $\kappa_1(G)=1$, by Theorem \ref{MT2-dl1}, $\Omega_1(G')$ is minimal normal in $G$.
It follows by Lemma \ref{CG(a)=1} that $C_{G'}(x)=1$ and so $C_{P}(x)=1$.
This implies that $P=[P,\lg x\rg]$ and so $$P=[P,\lg x\rg]\leq [N,H]\leq H\cap G'=H'.$$
Thus $N\cap G'=P=H'$.

Moreover, since $H\ntrianglelefteq G$, we have $P=H'<G'$ and so $G'\nleq N$. Since $G$ is 2-transitive on $\A_1(G)$, $G$ is primitive on $\A_1(G)$.
By \cite[Kapitel, Satz 1.4]{Hup}, $N$ is maximal in $G$. Thus $G=G'N$.
Notice that $H'=P=N\cap G'\unlhd N$ and $G'$ is abelian. We get $H'\unlhd G$.
Since $\kappa_1(G)=1$, by Theorem \ref{MT2-dl1}, $\Omega_1(G')$ is minimal normal in $G$. It follows that $H'=\Omega_1(G')$.

Now, since $\Omega_1(G')=H'<G'$, $\exp(G')\geqslant p^2$.
On the other hand, assume $N\leq H'\rtimes Q^g$ for some $g\in G$. Then $N<\Omega_2(G')\rtimes Q^g$.
Since $N$ is maximal in $G$, $G=\Omega_2(G')\rtimes Q^g$ and hence $G'\leq \Omega_2(G')$. That means $\exp(G')\leqslant p^2$. Thus $\exp(G')=p^2$.

\medskip
{\bf Step 5.} $q=2$ and $p$ is a Fermat prime, or $p=2$ and $q$ is a Mersenne prime.

\smallskip
Let $\Omega=\A_1(G)\setminus \{H\}$.
Then, by \cite[Kapitel II, Hilfssatz 1.8]{Hup}, $N$ is transitive on $\Omega$. By  Step 4, $G'\nleq N$.
Take $g\in G'\setminus N$.  Then $H^g\neq H$. It follows that $H^g\in\Omega$.
Since $G'$ is abelian,  $P\leq H^g$. Hence $P\in \Syl_p(H^g)$.

Let $N_1=N_N(H^g)$.  Then $|\Omega|=|N:N_1|$. Since $\kappa_1(G)=1$, $|\A_1(G)|=|G:N|$. Thus $$|G:N|=|N:N_1|+1.$$
Let $|G:N|=p^uq^v$. By Step 4,
$G' \cap N=\O_1(G')$. By Theorem \ref{MT2-dl1}, $G'$ is homocyclic. It follows that $p^u=|G':\O_1(G')|=|\O_1(G')|$.
Notice that $P\in{\rm Syl}_p(N)$ and $P\leq H^g$. We have $p\nmid |N:N_1|$.
Let $|N:N_1|=q^w$.
Then $p^uq^v=q^w+1$. It follows that $v=0$.
If $u=1$, then $p=q^w+1$. Thus $q=2$ and $p$ is a Fermat prime. If $w=1$, then $q=p^u-1$. Thus $p=2$ and $q$ is a Mersenne prime.

Assume $u,w\geqslant 2$. Then Catalan equation $x^u-y^w=1$ has no other solution in positive integers except $3^2-2^3=1$  by \cite[Theorem 5]{Miha}. Thus $p=3,u=2$, $q=2$ and $w=3$.
Then $|\O_1(G')|=3^2$.
Take $H=\O_1(G')\rtimes\lg g\rg\in\A_1(G)$, where $g\in Q\backslash Z(G)$. Then $g$ induces an automorphism of $\O_1(G')$ of order 2.
Let $a\in \O_1(G')$ and $a\ne 1$. If $a^g\neq a^{-1}$, then $aa^g\neq 1$. Clearly, $(aa^g)^g=aa^g$. That means $C_G'(g)\ne 1$. This contradicts
Lemma~\ref{CG(a)=1}.
If $a^g=a^{-1}$, then $\lg a,g\rg\in\A_1(G)$. But $|\lg a,g\rg'|=3<|\O_1(G')|$. This contradicts Step 4.
This means the case $u,w\geqslant 2$ doesn't occur.  \qed


\medskip
In the end of this section we give an example to explain that there exist the groups satisfying the conditions in Theorem \ref{MT2-dl2}.

\begin{exam} For any Fermat prime $p=2^{2^n}+1$, let $G=\lg a\rg\rtimes \lg b \rg$, where $\lg a\rg\cong\C_{p^2}$ and $b\in\Aut(\lg a\rg)$ such that $o(b)=p-1$.
For any Mersenne prime $q=2^{r}-1$, let $G=P\rtimes \lg b \rg$, where $P\cong\C_{2^2}^r$, $\lg b\rg\cong\C_{q}$ and $b$ acts irreducibly on $\Omega_1(P)$. In either case, we have $\kappa_1(G)=1$, $G$ is $2$-transitive on $\A_1(G)$ and $G$ is a Frobenius group.
\end{exam}


\section{Finite $p$-groups $G$ with $\kappa_1(G)=p$}

For a positive integer $t$, a finite non-abelian $p$-group $G$ is called an $\mathcal{A}_t$-group if its subgroups of index $p^t$ are abelian and it has at least one non-abelian subgroup of index $p^{t-1}$, which is introduced by Berkovich in \cite{BE}. Obviously, an $\mathcal{A}_1$-group is exactly a minimal non-abelian $p$-group, and a finite non-abelian $p$-group $G$ can be regarded as an $\mathcal{A}_t$-group for some positive integer $t$. We determine the finite non-abelian $p$-groups $G$ with $\kappa_{1}(G)=p$ in terms of $\mathcal{A}_t$-groups. This solves a problem posed by Berkovich, see \cite[Problem 23]{Be1}. For convenience, we introduce the following notations.

$\beta_1(G)$: the number of non-abelian subgroups of index $p$;

$Conj(G,H)=\{H^g\mid g\in G\}$: the set of conjugate class of subgroup $H$ in $G$;

$M\lessdot G$: $M$ is a maximal subgroup of $G$.

\begin{lemma}\label{beta1}
 Assume $G$ is a finite $p$-group which is neither abelian nor minimal non-abelian, that is, $G$ is an $\mathcal{A}_t$-group with $t\geqslant 2$. Let
$${\cal CA}_1(G)=\{Conj(G,H)\di H\in\A_1(G)\},\ N\Gamma_1(G)=\{M\lessdot G\di M'\neq 1\}$$
and $$\phi: Conj(G,H)\mapsto M, \ \forall\ H\leq M\lessdot G\ and\ H\in{\cal A}_1(G).$$
If the intersection of any two distinct maximal subgroups of $G$ is abelian,
then $\kappa_1(G)\geqslant\beta_1(G)$ and the equation holds if and only if $\phi$ is a bijection.
\end{lemma}

\demo
Since the intersection of any two distinct maximal subgroups of $G$ is abelian, it is a routine matter to get
$\phi$ is a surjection from ${\cal CA}_1(G)$ to $N\Gamma_1(G)$. Thus $\kappa_1(G)\geqslant\beta_1(G)$ and the equation holds if and only if $\phi$ is a bijection, 
\qed

\begin{lemma} \label{jdzq1}
Let $G$ be a finite $p$-group, $H\leq M\lessdot G$.
If $Conj(G,H)\neq Conj(M,H)$, then $$|Conj(G,H)|=p\cdot|Conj(M,H)|$$ and $Conj(G,H)$ is a partition of $p$ conjugate classes in $M$.
\end{lemma}
\demo It suffices to show that $N_G(H)=N_M(H)$. Since $H\leq M\lessdot G$, $$|N_G(H):N_M(H)|=|N_G(H):(N_G(H)\cap M)|\leqslant p.$$
If $|N_G(H):N_M(H)|=p$, then $$|Conj(M,H)|=|M:N_M(H)|=|G:N_G(H)|=|Conj(G,H)|.$$
It follows that $Conj(G,H)=Conj(M,H)$. This contradicts the hypothesis. \qed

\begin{lemma} \label{jdzq2}
Let $G$ be a finite $p$-group such that the intersection of any two distinct maximal subgroups of $G$ is abelian.
If $\kappa_1(G)=\beta_1(G)$, then $\kappa_1(M)\leqslant p$ for any maximal subgroup $M$ of $G$.
\end{lemma}
\demo Since $\kappa_1(G)=\beta_1(G)$, by Lemma \ref{beta1}, we get
$$\phi: Conj(G,H)\mapsto M,\ \forall\ H\leq M\lessdot G\ and\ H\in{\cal A}_1(G)$$
is a bijection from ${\cal CA}_1(G)$ to $N\Gamma_1(G)$. Thus ${\cal A}_1(M)=Conj(G,H)$.
It follows by Lemma \ref{jdzq1} that  $\kappa_1(M)=1$ or $p$ for any non-abelian $M\lessdot G$. Thus the result holds. \qed


\begin{cor}\label{zyc}
Let $G$ be a finite $p$-group. If $\kappa_1(G)\leqslant p$, then $\kappa_1(M)\leqslant p$ for any $M\lessdot G$.
\end{cor}

\demo By Lemma \ref{geqp}(2), $\kappa_1(G)=0, 1$ or $p$. 
If $\kappa_1(G)=1$, then, by Lemma \ref{geqp}(1), $G$ is a minimal non-abelian $p$-group and so $\kappa_1(M)=0$.
If $\kappa_1(G)=p$, by \cite[Theorem 17]{BE}, $d(G)=2$ and $G$ has an abelian subgroup of index $p$. 
Thus $\beta_1(G)=p$ and the intersection of any two distinct maximal subgroups of $G$ is abelian. By Lemma \ref{jdzq2}, $\kappa_1(M)\leqslant p$. \qed

%
%
%
%
%
%
%

\begin{lemma}{\rm(\cite[$\S 1$, Exercise 4]{BE1})}\label{jdl2}
Let $G$ be a finite $p$-group of maximal class and $G$ has an abelian subgroup $A$ of index $p$. Then all non-abelian subgroups of $G$ are of maximal class.
\end{lemma}

\begin{lemma}\label{jdlNGH}
Let $G$ be a finite $p$-group of maximal class with an abelian subgroup $A$ of index $p$ and $H<G$. If $H\nleq A$, then $|N_G(H):H|=p$.
\end{lemma}

\demo If not, then $|N_G(H):H|\geqslant p^2$. Take a subgroup $L$ of $G$ such that $$H\leq L\leq N_G(H)\ {\rm and}\  |L:H|=p^2.$$ Then $H\unlhd L$.
Since $H\nleq A$, we have $L\nleq A$ and $|L|\geqslant p^3$. It follows that $G=AL$ and $|L\cap A|\geqslant p^2$.
If $L$ is abelian, then $L\cap A\leq Z(G)$ and so $|Z(G)|\geqslant p^2$.
This contradicts that $G$ is of maximal class. Thus $L$ is not abelian.
By Lemma \ref{jdl2}, $L$ is of maximal class. It follows that $H=L'$ and so $H\leq G'\leq A$. This contradicts the hypothesis $H\nleq A$.  \qed

\begin{lemma}{\rm\cite[Lemma 3.2]{XAZ}}\label{Dp2-Lem1}
Let $G$ be a non-abelian two-generator $p$-group having an abelian maximal subgroup $A$. Then $G/Z(G)$ is of maximal class and $Z(M)=Z(G)$ for any non-abelian maximal subgroup $M$ of $G$.
\end{lemma}

\begin{theorem} \label{Main}
Assume $G$ is a finite $p$-group which is neither abelian nor minimal non-abelian, that is, $G$ is an $\mathcal{A}_t$-group with $t\geqslant 2$.
Then $\kappa_{1}(G)=p$ if and only if $G$ has an abelian subgroup of index $p$ and all non-abelian subgroups of $G$ are generated by two elements.
\end{theorem}

\demo $\Longrightarrow$: By \cite[Theorem 17]{BE}, $d(G)=2$ and $G$ has an abelian subgroup of index $p$.
Let $H$ be a non-abelian subgroup of $G$.
Then $\kappa_{1}(H)\leqslant p$ by Corollary \ref{zyc}. Thus $\kappa_{1}(H)=1$ or $p$ by Lemma \ref{geqp}.
If $\kappa_{1}(H)=1$, then $H$ is a minimal non-abelian $p$-group by Lemma \ref{geqp}(1) and so $d(H)=2$ by Lemma \ref{minimal non-abelian equivalent conditions}.
If $\kappa_{1}(H)=p$, then $d(H)=2$ by \cite[Theorem 17]{BE} again.

\medskip
$\Longleftarrow$: Assume $t=2$. Then all minimal non-abelian subgroups of $G$ are of index $p$ and so they are normal in $G$. It follows that $\kappa_{1}(G)=\alpha_1(G)=\beta_1(G)$.
Since $d(G)=2$ and $G$ has an abelian subgroup of index $p$, by Lemma \ref{Miller}, $\beta_1(G)=p$. So $\kappa_{1}(G)=p$.

\medskip
Assume $t>2$ and $K$ is an $\A_1$-subgroup of $G$. Let $H$ be a non-abelian subgroup of $G$. Then $d(H)=2$ and $H\cap A$ is an abelian maximal subgroup of $H$. 
This means that the hypothesis $``\Longleftarrow"$ in Theorem \ref{Main} is inherited by its non-abelian subgroups.
Following, we complete the proof by three steps.

\medskip
{\bf Step 1.} $Z(K)=Z(G)$ and $|G:K|=p^{t-1}$.

\smallskip

Since $t>2$, there exists non-abelian maximal subgroup $M$ of $G$ such that $K<M$. 
By inductions on $|G|$, we get $Z(K)=Z(M)$. By Lemma \ref{Dp2-Lem1}, $Z(M)=Z(G)$ and so $Z(K)=Z(G)$. 
Since $K$ is an $\A_1$-subgroup of $G$, by Lemma \ref{minimal non-abelian equivalent conditions}, $|K:Z(K)|=p^2$.
It follows that $|K:Z(G)|=p^2$. Thus all $\A_1$-subgroups of $G$ are of same order.
By \cite[Lemma 4.3]{Zjq}, $|G:K|=p^{t-1}$.


%

\medskip
{\bf Step 2.}  $\a_1(G)=p^{t-1}$.
\smallskip

By Step 1, all $\A_1$-subgroups are of index $p^{t-1}$ in $G$. 
It follows by \cite[Lemma 4.3]{Zjq} that $M$ is an $\mathcal{A}_{t-1}$-group for any non-abelian subgroup $M\lessdot G$. By inductions on $t$, $\a_1(M)=p^{t-2}$.
Since $d(G)=2$ and $G$ has an abelian subgroup, by Lemma \ref{Miller}, $G$ has $p$ non-abelian maximal subgroups. 
By the enumeration principle of P. Hall (see \cite[Theorem 5.2]{BE1}), $\a_1(G)=p\cdot\a_1(M)=p^{t-1}.$

\medskip
{\bf Step 3.}  $\kappa_1(G)=p$.
\smallskip

Let $\overline{G}=G/Z(G)$. Then $\overline{G}$ is of maximal class by Lemma \ref{Dp2-Lem1}. Notice that $K$ is not abelian. We get $K\nleq A$. By Step 1, $Z(K)=Z(G)$ and so $\overline{K}\nleq \overline{A}$.
It follows by Lemma \ref{jdlNGH} that $|N_{\overline{G}}(\overline{K}):\overline{K}|=p$ and so $|N_G(K):K|=p$. Thus
$$|Conj(K,G)|=|G:N_G(K)|=\frac{1}{p}|G:K|=p^{t-2}.$$
By the arbitrariness of $K$, we get $\kappa_1(G)=\frac{\a_1(G)}{|Conj(K,G)|}$. By Step 2, $\a_1(G)=p^{t-1}$ and so $$\kappa_1(G)=\frac{\a_1(G)}{|Conj(K,G)|}=\frac{p^{t-1}}{p^{t-2}}=p.$$
The proof is completed. \qed

\medskip
Finite $p$-groups $G$ whose non-abelian proper subgroups are generated by two elements have been determined in \cite{XAZ}.
If $G$ has an abelian subgroup of index $p$, then, by a trivial checking,  $\kappa_{1}(G)=p$ if and only if $G$ is one of the groups listed in \cite[Theorem 3.12 and 3.13]{XAZ} with the assumption that $m\geqslant 1$ instead of that $m\geqslant 2$.


\medskip

{\bf Acknowledgments.} The authors thank Professor Qinhai Zhang for his valuable suggestion on this theme.

\end {document}